\documentclass{article}

\usepackage[small,nohug,heads=vee]{diagrams}
\diagramstyle[labelstyle=\scriptstyle]
\usepackage{amsmath,amscd}

\begin{document}

Tsemo Aristide

College Boreal, 1 Yonge Street, 

M5E 1E5, Toronto, ON 

Canada

tsemo58@yahoo.ca

\bigskip
\bigskip

{\bf Closed models, strongly connected components and Euler graphs.}

\bigskip
\bigskip

\centerline {\bf Abstract.}

\bigskip

{\it In this paper, we  continue our study of closed models defined in categories of graphs. We  construct a closed model defined in the category of directed graphs which characterizes the strongly connected components. This last notion  has many applications,  and it plays an important role in the web search algorithm of Brin and Page, the foundation of the search engine Google. We also show that for this closed model, Euler graphs are particular examples of cofibrant objects. This enables us to interpret in this setting the classical result of Euler which states that a directed graph is Euleurian if and only if the in degree and the out degree of every of its nodes are equal. We also provide a cohomological proof of this last result.}

\bigskip
\bigskip

{\bf 1. Introduction.}

\medskip

In this paper, we pursue our investigation of closed models defined in the category $Gph$ of directed graphs. Recall that in [2] and [3], that we have published in collaboration with Terrence Bisson, we have introduced two closed models: the first is related to the zeta function of directed graphs and the second to dynamical systems. These constructions have been generalized in [10] where we have defined the notion of closed models defined by counting and study the existence of such closed models in the category of undirected graphs.  For the closed model defined in [2], a morphism of $Gph$ $f:X\rightarrow Y$ is a weak equivalence if and only if for every cycle $c_n, n>0$, the  morphism of sets $Hom(c_n,X)\rightarrow Hom(c_n,Y)$ induced by $f$ is a bijection.

In this paper,  we modify this condition by allowing $n$ to be equal to zero, otherwise said, we are counting also the nodes of $X$. This new closed model defined in $Gph$ enables to study other interesting properties of this category in particular it enlightens the important notion of strongly connected component of a directed graph, which has many applications in web search engines: the well known search engine Google designed by Brin and Page [5] uses the notion of pagerank to construct an hierarchy of the web which can be calculated by using strongly connected components and Markov matrices. More precisely, we show that a morphism $f:X\rightarrow Y$ is a weak equivalence for this closed  model if and only if it induces a bijection between the respective sets of strongly connected components of $X$ and $Y$ and its restriction to each strongly connected component of $X$ is an isomorphism onto a strongly connected component of $Y$. 
The cofibrant objects obtained here enable us also to study Eulerian graphs and to interpret the famous Euler theorem which states that a finite directed graph $X$ is Eulerian if and only if for every node $x$ of $X$ the inner and the outer degree of $x$ are equal. We also provide a construction of new closed models from a closed model defined by counting. This enables us to give a conceptual formulation of the closed model defined in [2].

We also introduce an homology theory in the category $Gph$ and show that the positive cycles of the first homology group of a directed graph is the set of cycles; this also enables us to give an homological interpretation of the Euler's theorem that we have just quoted
and to establish a link between the notions studied in this paper and simplicial sets. In this regard, we show that there exists a closed model defined in the category of $1$-simplicial sets also called the category of reflexive graphs which has many similarities which the closed model studied earlier in this paper.

\bigskip

{\bf 2. Some basic properties of the category of directed graphs.}

\medskip

Let $C$ be the category which has two objects that we denote by $0$ and $1$; the morphisms of $C$ which are not identities
are $s,t\in Hom(0,1)$.

\medskip

{\bf Definitions 2.1.}

The category $Gph$ of presheaves over $C$ is the category of directed graphs. Thus, a directed graph $X$
is defined by two sets $X(0)$ and $X(1)$, and two maps $X(s),X(t):X(1)\rightarrow X(0)$. The elements of $X(0)$ are
called the nodes of $X$ and the elements of $X(1)$ the arcs of $X$. For every arc $a\in X(1)$, $X(s)(a)$ is the source
of $a$ and $X(t)(a)$ is the target of $a$. We will also often say that $a$ is an arc between $X(s)(a)$ and $X(t)(a)$ or that $a$
connects $X(s)(a)$ and $X(t)(a)$.

A morphism $f:X\rightarrow Y$ between two directed graphs is a morphism of presheaves: it is  defined by two maps
$f(0):X(0)\rightarrow Y(0)$ and $f(1):X(1)\rightarrow Y(1)$ such that $Y(s)\circ f(1) = f(0)\circ X(s)$ and $Y(t)\circ f(1) = f(0)\circ X(t)$.

Let $X$ be a finite directed graph, suppose that the cardinality of $X(0)$ is $n$, the adjacency matrix $A_X$ of $X$ is the $n\times n$ matrix whose entry $(i,j)$ is the cardinal of $X(x_i,x_j)$, the set of arcs between $x_i$ and $x_j$.

\medskip

{\bf Definitions 2.2.}

Let $X$ be a graph, and $x$ a node of $X$. We denote by $X(x,*)$ the set of arcs of $X$ whose source is $x$, and by $X(*,x)$
the set of arcs of $X$ whose target is $x$. If $X$ is finite, the inner degree of $x$ is the cardinality of $X(*,x)$ and the outer degree
of $x$ is the cardinality of $X(x,*)$.

Examples of directed graphs are:

The directed dot graph $D$; $D(0)$ is a singleton and $D(1)$ is empty. Geometrically it is represented by a point.

The directed arc $A$. The set of nodes of $A$ contains two elements $x,y$, and $A$ has a unique arc $a$ such $A(s)(a)=x$
and $A(t)(a)=y$. Geometrically, it is represented by an arc between $x$ and $y$ as follows: $x \longrightarrow y$.

The directed cycle $c_n, n\geq 1$ of length $n$; $c_n(0)$ is a set which contains $n$ elements that we denote by $x^n_0,...,x^n_{n-1}$. For $i<n-1$, there is
a unique arc $a^n_i$ whose source is $x^n_i$ and whose target is $x^n_{i+1}$; there is an arc $a^n_{n-1}$ whose source is $x^n_{n-1}$ and whose target
is $x^n_0$. Often, we will say that $D$ is the cycle $c_0$ of length $0$.

\medskip

The directed line $L$ is the graph such that $L(0)$ is the set of  integers $Z$, and for every integer $n$,
there exists a unique arc $a_n$ such that $L(s)(a_n) = n$ and $L(t)(a_n) = n+1$.

\medskip

The directed path $P_n$ of length $n$; the set of nodes $P_n(0)$ has $n$ elements $x^n_0,....,x^n_{n-1}$ and for $i<n-1$, there
exists an arc $a^n_i$ between $x^n_i$ and $x^n_{i+1}$; $x^n_0$ is the source of the path and $x^n_{n-1}$ is its end.

\medskip

{\bf Definitions 2.3.}

Let $X$ be an object of $Gph$ and $x,y$ two nodes of $X$. A path between $x$ and $y$ is  a morphism $f:P_n\rightarrow X$ such that $f(0)(x^n_0) = x$ and $f(0)(x^n_{n-1}) = y$. We say that $X$ is connected if and only if for every nodes $x$ and $y$ of $X$, there exists a finite set of nodes $(x_i)_{i=1,...,l}$
such that $x_1=x, x_l=y$ and for $i<l$, there exists a path between $x_i$ and $x_{i+1}$ or a path between $x_{i+1}$ and $x_i$.

The graph $X$ is strongly connected if and only if for every nodes $x,y$ of $X$ there exists a path between $x$ and $y$ and a path between
$y$ and $x$. This is equivalent to saying that there exists a cycle which contains $x$ and $y$.

Let $X$ be a directed graph, consider the equivalent relation $R$ defined on the space of nodes of $X$ such that $x Rx$
for every $x\in X(0)$, if $x$ is distinct of $y$ then $x R y$ if and only if  there exists a cycle which contains $x$ and $y$.

We denote by $U_1,...,U_p,...$ the set of equivalent classes of this relation.  We denote by $X_{U_i}$ the subgraph of $X$ whose set
of nodes is $U_i$. An arc $a\in X(1)$ is an arc of $X_{U_i}$ if and only if $X(s)(a)$ and $X(t)(a)$ are elements of $U_i$.
The graphs $X_{U_i}$ are the strongly connected components of $X$.

\bigskip

{\bf 3. Closed models in $Gph$.}

\medskip

We recall now the notion of closed model category:

\medskip

{\bf Definition 3.1.}

Let $C$ be a category, we say that the morphism $f:X\rightarrow Y$ has the left lifting property with respect to   the morphism $g:A\rightarrow B$ (resp., $g$ has the  right lifting property with respect to $f$) if and only if for each commutative square

$$\dot{}
\begin{CD}
X  @> l>> A\\
@VV f V    @VV g V\\
Y @>m>>  B
\end{CD}
$$
there exists a morphism $n:Y\rightarrow A$ such that $ l = n\circ f$ and $m = g\circ n$.
Let $I$ be a class of maps of $C$, we denote by $inj(I)$ the class of morphisms of $C$ such that for every $f$ in $I$ and
every $g\in inj(I)$, $g$ has the right lifting property with respect to $f$.
We denote $cell(I)$ the subclass of maps of $C$ which are retracts of transfinite composition of pushouts of elements $I$.

Two class of maps $L$ and $R$ define a weak factorization system $(L,R)$ of $C$ if and only if:
 for every morphism $f$ of $C$, there exists $g\in R$
and $h\in L$ such that $f = g\circ h$ and  $L$ is the class of morphisms which have the left lifting property with respect to every morphism $R$ and
$R$ is the class of morphisms which have the right lifting property with respect to every morphism of $L$.

\medskip

{\bf Definition 3.2.}

A closed model category is a category $M$ which has projective limits and inductive limits
endowed with three subclasses of morphisms $W,F,C$ called respectively the weak equivalences, the fibrations
and the cofibrations. We denote by $F'$ (resp., $C'$) the intersection $F\cap W$ (resp., $C\cap W$). The subclass
$F'$ is called the class of weak fibrations and $C'$ the class of weak cofibrations.  The following
 two axioms are also satisfied:

M1. $(C,F')$ and $(C',F)$ are weak factorization systems.

M2 Let $f:X\rightarrow Y$ and $g:Y\rightarrow Z$ be two maps in $M$, if two maps of the triple $\{ f,g,g\circ f\}$ is a weak equivalence so is the third.

In this paper, we are only going to consider locally presentable categories. This has the virtue to avoid set theoretical difficulties
when one tries to find weak factorizations systems. We are going to use Proposition 1.3 of Beke [1] which asserts that if $I$ is a class of morphisms
of a locally presentable category, $(cell(I), inj(I))$ is a weak factorization system. The categories of graphs used here are
locally presentable categories since they are isomorphic to categories of presheaves defined on a small category.

\medskip

{\bf Definition 3.3.}

A  closed model structure defined on $C$ is cofibrantly generated if and only if there exists a set of morphisms $I$ (resp., $J$)  such that
 $inj(I)$ (resp., $inj(J)$) is the class of weak fibrations (resp., the class of fibrations).

\medskip

In [10] we have introduced the notion of a closed model category defined by counting which we outline: it is a closed model category $C$,   whose class of weak equivalences
$W$ is defined as follows:

 Firstly, we consider a set of objects of $C$, $(X_l)_{l\in L}$. Let $\phi$ be the initial object of $C$, we can define the morphisms
 $i_l:\phi\rightarrow X_l$ and the folding morphism $j_l:X_l+X_l\rightarrow X_l$ which is the sum of two copies of $Id_{X_l}:X_l\rightarrow X_l$. The class $W$ is $inj(I)$ where  $I=\{i_j,j_l; {l\in L}\}$. Thus a morphism $f:X\rightarrow Y$ is a weak equivalence
  if and only if for every $l\in L$, the map $Hom(X_l,X)\rightarrow Hom(X_l,Y)$ which sends $g:X_l\rightarrow X$ to $f\circ g$ is bijective, and $(cell(I),W)$ is a weak factorization system.
 We can define a closed model on $C$ whose class of weak equivalences is $W$, the class of fibrations is the class of morphisms of $C$
 and the class of cofibrations is $cell(I)$. Remark that such a closed model is cofibrantly generated since its class of fibrations is $inj(\phi)$ where $\phi$ is the initial object.

 \medskip

 {\bf Proposition 3.1.}

 {\it Let $W$ be the class of weak equivalences of the closed model defined by counting the objects $(X_l)_{l\in L}$,
  $J$  a set of morphisms $(f_j)_{j\in P}$
  such that $cell(J) \subset W$.  Denote by $F$ the class $inj(J)$
 and by $Cof$ the class of morphisms $cell(I\bigcup J)$. Then $(W,F,Cof)$ defines a closed cofibrantly generated closed model on $C$.}

  \medskip

  {\bf Proof.}

   We are going to apply the result of D. Kan quoted by Hirschhorn [9] p. 213, Theorem 11.3.1 that shows that the sets of morphisms $I\bigcup J$ and $J$ define a cofibrantly
   generated closed model on $C$ where $I=\{i_l,j_l; {l\in L}\}$.

   A morphism $f$ of $cell(J)$ is an element of $W$ by assumption, and is obviously contained in $cell(I\bigcup J)$. A morphism
    $f$ of $C$ which is right orthogonal to $I\bigcup J$ is a weak equivalence since it is right orthogonal to $I$ and is obviously right orthogonal
    to $J$. This verifies the conditions 2 and 3 of the theorem of Kan. A morphism $f$ which is right orthogonal to $J$ and is in $W$ is a morphism
    right orthogonal to $I\bigcup J$. This verifies the condition $4 (b)$.

\medskip

 {\bf Examples.}

 We present now the following closed model defined  by counting the cycle graphs $(c_n)_{n>0}$ in the category $Gph$. A morphism $f:X\rightarrow Y$
 is contained in the class $W'$ of    weak
 equivalences of this closed model if for every $n>0$, the map $Hom(c_n,X)\rightarrow Hom(c_n,Y)$
 is bijective. We have a closed model $(W',Fib',Cof')$ for which $Fib'$ is the class of all the maps and $Cof'$ is
 $cell(i_n,j_n,n>0)$, where $i_n:\phi\rightarrow c_n$ and $j_n:c_n+c_n\rightarrow c_n$.

 We can apply the Proposition 3.1, to obtain other closed models with the same class of weak equivalences. On this purpose,
 consider a non empty graph $X$ such that for every integer $n>0$, $Hom(c_n,X)$ is empty. Such a graph is called acyclic. Let $x$ be any node of $X$,
 consider the morphism $s^x:D\rightarrow X$ such that the image of $s^x(0)$ is $x$. An element of $cell(s^x)$ is a composition of morphisms  $f:Y\rightarrow Z$, where $f$ is the canonical embedding of $Y$ into a graph  $Z$  obtained by attaching an acyclic graph to a node of $Y$. See also [4] Proposition 4. We deduce that the class $cell(s^x)$ is contained in $W'$. We
  can thus apply the Proposition 3.1 to obtain the closed model $(W',F_X,Cof_X)$ such that $F_X$ is $inj(s^x)$, and
  $Cof_X = cell(i_n,j_n,s^x,n>0)$. In particular, if $s:D\rightarrow A$ is the morphism between the dot graph and the arc
  graph such that $s(0)$ is the source of $A$, we obtain the closed model presented in [2] for which the class of fibrations is $inj(s)$
  and the cofibrations are $cell(i_n,j_n,s)$.
  Other examples may rise some interest. We can define $t:D\rightarrow A$ such that the image of $t(0)$ is the target of $A$ and obtain
  a closed model whose weak equivalences are $W'$, the class of fibrations is $inj(t)$ and the cofibrations are $cell(i_n,j_n,t,n>0)$.
  We can also defined the closed model whose weak equivalences are $W'$, the class of fibrations is $inj(s,t)$
  and the cofibrations are $cell(i_n,j_n,s,t)$.

 \bigskip

 {\bf 4. Closed models and strongly connected components.}

 \medskip

 One of the main purposes of this paper is to study a closed model defined by counting on $Gph$ related to $(W',Fib',Cof')$. This time,  we count the cycles $(c_n)_{n\geq 0}$. That
 is, we are also counting  nodes.
 Thus a weak equivalence $W$ for this closed model is a morphism $f:X\rightarrow Y$ such that for every $n\geq 0$, the map $Hom(c_n,X)\rightarrow Hom(c_n,Y)$ which associates $f\circ g$ to each element $g\in Hom(c_n,X)$ 
 is bijective.  We  obtain a closed model $(W,Fib,Cof)$ for which $Fib$ is the class of all the morphisms of $Gph$ and $Cof$ is
 $cell(i_n,j_n, n\geq 0)$. This closed model is related to strongly connected components of directed graphs, a notion which is intensively used in computer science and in particular in web search as shows the work of Brin and Page [5], the conceptual foundation of the search engine Google. Given a network (a directed graph), it is important for a web search engine to recommend pages to an user, on this purpose, a weight is assigned to each page (vertex) called the pagerank which depends on the number of important links that the page receives (the weight of the source of the incoming arcs). If $A$ is the adjacency matrix of the network, to obtain the pagerank, one has to define a new matrice $P$ by replacing the non zero coefficients of $A$ by numbers which quantify the importance of the link, and the pagerank of the page $i$ is just the sum of the entries of the $i$-row of  $P$. It is also reasonable to normalize the columns of the matrix $P$ to minimize the importance of outgoing links from a page, so surfing online is assimilated to a random walk described by the Markov matrix $P$. Linear algebra shows thus the pagerank is an eigenvalue of $P$. If $P$ is irreducible, the Perron theorem shows the existence of a  unique maximal positive eigenvalue which defines the pagerank. The fact that $P$ is irreducible means also that the graph is strongly connected. In practice this is not true, but research shows that 90 percent of the world wide web is connected and contains a giant strongly connected component. To cope of the general situation, google uses transition probabilities. 
 
 We have  the following result:

\medskip

{\bf Theorem 4.1.}

{\it A morphism $f:X\rightarrow Y$ of $Gph$ is an element of $W$ if and only if it induces a bijection between the sets of  strongly connected components
of $X$ and $Y$ and the restriction of $f$ to a strongly component of $X$ is an isomorphism onto a strongly connected component of $Y$.}

\medskip

{\bf Proof.}

 Firstly, we show that the image of a strongly connected component $U$ of $X$ is a strongly connected component. The restriction $f_{\mid U}$ of $f$ to $U$ is injective on nodes, since $f$ induces a
bijection on the set of nodes. Let $a$ and $b$ be two arcs of $U$ such that $f(1)(a) = f(1)(b)$. Since $f$ is injective on nodes,
$s(a) = s(b)$ and $t(a) = t(b)$. Consider a path $p$ in $U$ between $t(a)$ and $s(a)$. We can construct two cycles
$c$ and $c'$ obtained respectively by the concatenation of $a$ and $p$ and the concatenation of $b$ and $p$ The images of $c$ and
$c'$ by $f$ coincide. This implies that $c=c'$ since $f$ is injective on cycles, thus $a=b$. The image of $U$ is
thus imbedded in a strongly connected component $V$ of $Y$.

Suppose that there exists a node $y$ in $V$ which is not in the image of $U$. Let $y'=f(0)(x)$, $x\in U$. Since $V$ is strongly connected,
there exists a cycle $c$ of $V$ whose set of nodes contains $y$ and $y'$. Consider the cycle $c'$ of $X$ whose image by $f$ is $c$; $c'$
contains $x$ since $f(0)$ is injective. This implies that $c'$ is in $U$, and $c$ is contained in the image of $U$. This is a contradiction with the fact that $y$
is not in the image of $U$. Consider an arc $b$ of $V$ which is not in the image of $U$. There exists  a cycle $c$ of $V$ that contains
$b$. Since $f$ induces a bijection on  cycles, there exists a cycle $c'$ of $X$ whose image by $f$ is $c$.
Let $a$ be the arc of $c'$ whose image by $f$ is $b$; $s(a)$ and $t(a)$ are contained in $U$ since their image are contained in $f(0)(U)$.
This implies that $a$ is in $U$ since $U$ is a strongly connected component and henceforth $b$ is in the image of $U$. Thus the restriction of $f$ to $U$ is surjective on arcs. Since the restriction of $f$ to $U$ is injective, we deduce that $f$ induces an isomorphism of $U$ onto its image $V$.

Let $V$ be a strongly connected component of $Y$, and $y$ a node of $V$. There exists a node $x\in X(0)$ such that $f(0)(x) = y$. The
image of the strongly connected component which contains $x$ is $V$. This implies that $f$ induces a bijection on strongly connected components.

Conversely, suppose that $f$ induces a bijection between the set of on strongly connected components
of $X$ and $Y$ and the restriction of $f$ to a strongly connected component of $X$ is an isomorphism. Let $c$ and $c'$ two $n$-cycles ($n$ eventually $0$) of $X$ whose
image by $f$ coincide. This implies that that $c$ and $c'$ are in the same strongly connected component $U$, and are equal since the restriction
of $f$ to $U$ is an imbedding. Let $c$ be a cycle of $Y$, $c$ is an element of a strongly connected component $V$. The strongly
connected component $U$ of $X$ whose image maps isomorphically to $V$ contains a cycle whose image is $c$. We deduce that $f$ is a weak equivalence.

\medskip

{\bf Corollary 4.1.}

{\it A morphism $f:X\rightarrow Y$ between two strongly connected directed  graphs is a weak equivalence if and only if it is an isomorphism.}

\bigskip

{\bf Cofibrant replacement.}

\medskip

We are going to study in this section the notion of cofibrant replacement for the closed model defined in this section 4 on $Gph$ by $(W,Fib=Hom(Gph),Cof)$.
Recall that an object $X$ is cofibrant if and only if the map $\phi\rightarrow X$ is a cofibration where $\phi$ is the initial object. The object $Y$ is a cofibrant replacement of $X$ if and only if $Y$ is a cofibrant object and there exists a weak equivalence $f:Y\rightarrow X$. We know that $Cof = cell(i_n,j_n,n\geq 0)$.
This implies that the $n$-cycles $n\geq 0$ are cofibrant.  We deduce also that the sum of cycles are cofibrant objects.

\medskip

{\bf Some cofibrant maps: Gluing nodes and paths.}

\medskip

Let $f:c_0\rightarrow c_m$ and $g:c_0\rightarrow c_n$ two morphisms of graphs. Consider the  pushout diagram:

 $$\dot{}
\begin{CD}
c_0  @> f>> c_m\\
@VV g V    @VV  V\\
c_n @>m>>  X
\end{CD}
$$

The graph $X$ is   obtained by identifying a node of $c_m$ with a node of $c_n$. We say also that $X$ is obtained by attaching $c_m$ and $c_n$
by a node. The graph $X$ is cofibrant. We can iterate this operation to create more cofibrant objects: for example we can attach more cycles or identify
paths as follows:

Consider the graph $X$ defined as follows: there exist two cycles $c_m$ and $c_n$, nodes $x,y$ of $c_m$ and nodes $x',y'$ of $c_n$ such that
there exist  a path $p_1\in c_m$ between $y$ and $x$ and a path $p_2$ in $c_n$ between $y'$ and $x'$ which have the same length. We can construct the graph $X$ obtained by attaching $c_m$ and $c_n$ by identifying $x, x'$ and $y,y'$. We denote by
$[x]$ $(resp., [y])$ the node of $X$ corresponding to $x$ (resp., $y$). In $X$, we have  paths $l_1,l_2$ between $[y]$ and $[x]$ and obtained respectively from $p_1$ and $p_2$ and which have the same length. There exist also another $l_3$ between $[x]$
and $[y]$ in $X$. We can construct the cycles $c = l_1l_3$ and $l_2l_3$ which have the same length $p$.
Let $f:c_{p}\rightarrow X$ whose image is $l_1l_3$ and $g:c_{p}\rightarrow X$ whose image is $l_2l_3$. We can
construct the pushout of $f+g:c_p+c_p\rightarrow X$ by $j_{p}:c_{p}+c_{p}\rightarrow c_{p}$. It is a morphism $h:X\rightarrow Y$ and $Y$ is obtained from $X$ by identifying  $l_1$ and $l_2$. We say that $Y$ is obtained by gluing the paths $l_1$ and $l_2$.

\medskip

{\bf Theorem 4.2.}

{\it A strongly connected graph is a cofibrant object.}

\medskip

{\bf Proof.}

Let $X$ be a strongly connected graph. There exists a family of cycles $(c_{n_i},i\in I)$ and a morphism $f:\sum_ic_{n_i}\rightarrow X$
surjective on nodes and arcs. We can write $f = h\circ g$ where $g$ is a cofibration and $h$ a weak fibration. Write
$h:Y\rightarrow X$, without restricting the generality, we can suppose that the image $Y$ of $g$ is connected. Thus $Y$ can be constructed from a cycle $c_p$  by repeating the following operations: attach a cycle to a point, identifying two nodes or two arcs. This implies that $Y$  is strongly connected. The Corollary 4.1 implies that $h$ is an isomorphism, we deduce that $f$ is a cofibration and $X$ is cofibrant.

\medskip

The previous construction yields to the following:

\medskip

{\bf Corollary 4.2.}

{\it Let $X$ be a directed  graph, consider the subgraph $c(X)$ of $X$ which has the same nodes of $X$, an arc of $X$
is an arc of $c(X)$ if and only if it is contained in a strongly connected component of $X$, the canonical embedding $c_X:c(X)\rightarrow X$
is a cofibrant replacement of $X$.}

\medskip

{\bf Proof.}

The graph $c(X)$ is the disjoint union of the strongly connected components of $X$. The Theorem 4.2 implies that $c(X)$ is a cofibrant object, and the Theorem 4.1
implies that the canonical embedding $c(X)\rightarrow X$ is a weak equivalence.

\bigskip

{\bf Application to Eulerian graphs.}

\medskip

We are going to apply these results to Eulerian cycles. Remark that:

\medskip

{\bf Proposition 4.1.}

{\it Let $X$ be a finite strongly connected directed graph, there exists an integer $n(X)$, and a morphism $f:c_{n(X)}\rightarrow X$ surjective  on arcs.}

\medskip

{\bf Proof.}

We fix a node $x_0$ of $X$. We can index the arcs of $X$ by $a_1,...,a_l$. Since $X$ is strongly connected, there exists a path
$p_i$ from $x_0$ to $s(a_i)$ and a path $p_i'$ from $t(a_i)$ to $x_0$ $i=1,...l$. We can construct the cycle $p'_la_lp_l...p'_ia_ip_i...p'_1a_1p_1$
which contains all the arcs of $X$.

\medskip

This leads to to the following definition:

\medskip

{\bf Definition 4.1.}

An Eulerian cycle in a directed graph $X$ is a cycle $f:c_n\rightarrow X$ such that $f(1)$ is a bijection. 

\medskip

We have the following proposition:

\medskip

{\bf Proposition 4.2.}

{\it A finite directed graph $X$ is Eulerian if and only if it is cofibrant and obtained from a cycle  by identifying nodes.}

\medskip

{\bf Proof.}

Let $X$ be an Eulerian graph. There exist an integer $n$ and a  morphism $f:c_n\rightarrow X$ surjective on nodes and bijective on arcs; $f$ is a cofibration since it is the
composition of morphisms which identify nodes  and
henceforth, we deduce that $X$ is cofibrant since $c_n$ is cofibrant.  Conversely,
a cofibrant graph $X$ obtained from a cycle $c_n$ by identifying some of its nodes is Eulerian and the canonical morphism $f:c_n\rightarrow X$
is an Eulerian cycle.

\medskip

{\bf Proposition 4.3.}

{\it Consider a graph $X$ constructed recursively as follows:  $X_0$ is a cycle $c_n$, to construct $X_1$, identify two nodes of $c_n$ or attach  a cycle to a node of $c_n$. Suppose defined $X_n$, to obtain $X_{n+1}$,  identify two nodes of $X_n$ or attach a cycle to a node of $X_n$. Each
graph $X_n$ is Eulerian.}

\medskip

{\bf Proof.}

The graph $X_0=c_n$ is Eulerian. Suppose that $X_n$ is Eulerian. Let $f:c_p\rightarrow X_n$ be an Eulerian cycle. If $X_{n+1}$
is obtained from $X_n$ by identifying two nodes, let $g:X_n\rightarrow X_{n+1}$ be the identifying morphism, $g\circ f$ is an Eulerian
cycle of $X$. Suppose that $X_{n+1}$ is obtained from $X_n$ by attaching a cycle $c_m$. The concatenation of the cycles $f$ and $c_m$ is
an Eulerian cycle of $X_{n+1}$.

\medskip

{\bf Theorem 4.3.}

{\it A finite directed connected graph $X$ is obtained by the processus described in Proposition 4.3 if and only if for every node $x$ of $X$,   the in and out degree of $x$ are equal.}

\medskip

{\bf Proof.}

Suppose that $X$ is an Eulerian graph, then Proposition 4.2 shows that there exists a sequence of graphs $X_0=c_n,...,X_n=X$ such that $X_{i+1}$ is obtained from $X_i$ by identifying two nodes of $X_i$. The identification of two nodes of an Eulerian graph
 increases the in degree and the out degree of a node by the same number, we deduce that if  $X_i$ is Eulerian, then $X_{i+1}$ is Eulerian.  Since $c_n$ is Eulerian, we deduce recursively that $X$ is Eulerian. 

 Conversely, suppose that $X$ is a connected directed finite graph such that the in degree and the out degree of every node of $X$
 coincide, we are going to show that $X$ is constructed by the process described at Proposition 4.3. Let $x$ be any node of $X$ and $a_0\in X(x,*)$,
 then $X(t(a_0),*)$ is not empty since its in degree is equal to its out degree, we consider $a_1\in X(t(a_0),*)$, if $t(a_1) = x$ we stop otherwise
 there exists $a_2\in X(t(a_1),*)$ by continuing this process we obtain a cycle $f_1:c_{n_1}\rightarrow X$ injective on arcs. We can consider the subgraph
 $X_1$ of $X$ which is the image of $f_1$; $X_1$ is obtained from $c_{n_1}$ by identifying nodes. If $X_1$ is not $X$, since $X$ is connected,
 we have $x_2\in X_1$ such that $X(x_2,*)$ contains an arc $a^2_1$ which is not in $X_1$, since the in degree and the out degree of $t(a^2_1)$
 are equal, if $t(a^2_1)$ is distinct of $x_2$ there exists an arc $a^2_2\in X(t(a^2_2),*)$ as above, we conclude the existence
 of an injective morphism $f_2:c_{n_2}\rightarrow X$ whose image is a cycle through $x_2$. We can construct the subgraph of $X$ which is the union
 of $X_1$ and the image of $f_2$. Remark that $X_2$ is obtained from $X_1$ by attaching a cycle and identifying nodes.
 We can repeat the process to obtain an increasing sequence of graphs $X_1\subset X_2\subset...X_i\subset X_{i+1}\subset...$ such
 that $X_{i+1}$ is obtained from $X_i$ by attaching a cycle and identifying nodes of this cycle. Since $X$ is finite, we deduce the existence
 of $n$ such that $X_n= X$. The Theorem 4.3 shows that $X$ is Eulerian.

 \medskip

{\bf Corollary. 4.3. (Euler).}

{\it A finite directed graph $X$ is Eulerian if and only if for every node $x$ of $X$, the in and out degree of $x$ are equal.}

\medskip

{\bf 5. Cohomological interpretation.}

\medskip

Let $X$ be a directed graph. We denote by $Z(X(0))$ (resp., $Z(X(1))$ the free commutative group generated by the set $X(0)$ (resp., by the arcs of $X$).
The elements  of $Z(X(0))$ are called the $0$-chains.  A $1$-chain $u$ of $X$ is the linear sum $\sum_{i=1}^{i=l}d_if_{n_i}$ where $d_i$ is an integer and   $f_{n_i}:P_{n_i}\rightarrow X$ is a morphism between the path of length $n_i$ and $X$. We denote by $Z(ch(X))$ the space of $1$-chains of $X$. To each $1$-chain $u$, we associate $u'$ the element of $Z(X(1))$ defined by $\sum_{i=1}^{i=l} d_i\sum_{m=0}^{m=n_i-1}f_{n_i}(1)(a^{n_i}_m)$, we will often call $u'$ the image of $u$.
 
 We say that $u$ is positive if and only if $d_i\geq 0, i=1,...,l$.
 
 The length $l_X(u)$ of $u$ is $\sum_in_i\mid d_i\mid$.
 
  Suppose that $X$ is finite, for each arc $a\in X(1)$, we define the morphism
  $f_a:P_1\rightarrow X$ whose image is $a$; the fundamental chain $[X]$ of $X$ is $\sum_{a\in X(1)}f_a$.

We define the linear map $d^X_1:Z(ch(X))\rightarrow Z(X(0))$ such that for every chain $f:P_n\rightarrow X$ of $X$, $d^X_1(f) = t(f)-s(f)$. Remark that $d_1^X(f) =\sum_{i=0}^{i=n-1} t(f(1)(a^n_i))-s(f(1)(a^n_i))$.

 We also define the linear map $d^X_0:Z(X(0))\rightarrow Z$ such that for every node $x$ of $X$, $d^X_0(x) = 1$. We have the relation $d^X_0\circ d^X_1 = 0$. We denote
by $H_1(X)$ the kernel of $d_1$, and by $H_0(X)$ the quotient of the kernel of $d_0$ by the image of $d_1$.

Each morphism $f:X\rightarrow Y$ between directed graphs induces natural morphisms $f^*_0:Z(X(0))\rightarrow Z(Y(0))$ and $f_1^*:Z(ch(X))\rightarrow Z(ch(Y))$.

Remark that if $f:c_n\rightarrow X$ is an $n$-cycle of $X$, the composition
of $f\circ p_n$ of $f$ with the canonical morphism $p_n:P_{n+1}\rightarrow c_n$ is a chain such that $d^X_1(f\circ p_n)=0$.

\medskip

{\bf Proposition 5.1.}

{\it Let $X$ be a finite directed graph, $u=\sum_{i\in I} d_if_{n_i}$ a positive $1$-chain,
$d_1^X(u)=0$ if and only if there exists a finite set of cycles $g_j:c_{n_j}\rightarrow X$ such that the images of $\sum_i d_if_{n_i}$ and $\sum_j g_j\circ p_{n_j}$ coincide.}

\medskip

{\bf Proof.}

Without restricticting the generality, we can assume that $d_i=1, i\in I$ since the chain is positive. We are going to give a recursive proof depending of the cardinality of $I$.
Suppose that $I$ is a singleton, then $u=f$ where $f:P_n\rightarrow X$. The fact that $d_1^X(f)=0$ is equivalent to say that $f$ factors by a morphism $c_n\rightarrow X$.

Suppose that the result is true if the cardinality of $I$ is $l$. Assume now that the cardinality of $I$ is $l+1$. Remark that $d^X_1(u)=\sum_if_{n_i}(0)(t(P_{n_i}))-f_{n_i}(0)(s(P_{n_i})) = 0$. This implies the existence of $i_p$ such that
$f_{n_{i_p}}(0)(s(P_{n_{i_p}}))=f_{n_0}(0)(t(P_{n_0}))$ we can thus define the concantenation $f_{n_{i_p}}f_{n_0}$ of $f_{n_p}$  which is an $n_0+n_p$-chain.
We  consider
the family $L=\{f_{n_i}, f_{n_{i_p}}f_{n_0}, i\in I\}-\{f_{n_0},f_{n_{i_p}}\}$ whose cardinal is strictly inferior to the cardinal of $I$ and such that $\sum_{i\neq 0,p }f_{n_i}+f_{n_{i_p}}f_{n_0}$ has the same image than $u$. We can apply the recursive hypothesis to it and obtain a family of cycles $g_j:c_{n_j}\rightarrow X$ such that the images of $\sum_i d_if_{n_i}$ and $\sum_j g_j\circ p_{n_j}$ coincide. 

\medskip

This enables to give another proof of the theorem of Euler:

\medskip

{\bf Corollary. 5.1. (Euler).}

{\it Let $X$ be a finite connected directed graph, there exists a morphism $f:c_n\rightarrow X$ bijective on arcs if and only if for every node $x$ of $X$ the in and the out degrees of $x$ coincide.}

\medskip

{\bf Proof.}

Suppose that for every node $x$ of $X$, the in degree $in(x)$ and the out degree $out(x)$ of $X$ coincide, we have $d^X_1([X]) =\sum_{x\in X(0)}(out(x)-in(x)) =0$. The  Proposition 5.1 implies the existence of morphism $f_{n_1}:c_{n_1}\rightarrow X,..., f_{n_l}:c_{n_l}\rightarrow X$ such that $\sum_{i=1}^{i=l}f_{n_i}$ and $[X]$ have the same image. We also deduce that  $\sum_{i=1}^{i=l}f_{n_i}$ is bijective on arcs since the coefficients of its image are $1$. Since $X$ is connected, we deduce the existence of a morphism 
$f:c_n\rightarrow X$ bijective on arcs by making a concatenation of $f_{n_i},i=1,...l$.

{\bf Remark.}

Let $X$ be a finite graph, $H_0(X)=0$ if and only if $X$ is connected, and 
$H_1(X)=0$ if and only if $X$ is acyclic: this is equivalent to saying that
for every integer $n>0$, $Hom(c_n,X)$ is empty. In fact, there exists a bijection between the set of cycles of $X$ and positive elements of $H_1(X)$. This allows to give another description of the class of weak equivalences $W'$ studied in: a morphism $f:X\rightarrow Y$ is an element of $W'$ if and only if $f_1^*:H_1(X)\rightarrow H_1(Y)$ is bijective on positive chains.

 \medskip

  Let $X$ be a finite strongly connected finite directed graph. We have seen that there exists a morphism $f:c_n\rightarrow X$ surjective
  on nodes and arcs. A good question is to find the lower bound $n(X)$ of $n$. We know that if $X$ is Eulerian, $n(X)$ is the cardinal
  of the number of arcs of $X$. The Proposition 5.1 shows that to find $n(X)$, it is sufficient to find a positive chain $c$ such that $d_1([X]+c)=0$ and the length of $l([X]+c)$ is minimal.

\bigskip

{\bf 6. Closed models on $RGph$.}

\medskip

The cohomological interpretation of the proof ot the Euler theorem suggests that this theory is related to simplicial sets. In fact, $1$-simplicial sets are often called reflexive graphs, in this part, we are going to study a closed model in the category $RGph$ of reflexive graphs related to the closed model that we have just studied in $Gph$.
Consider the category $C_R$ which has two objects $0_R$ and $1_R$, the morphisms of $C_R$ different of the identities are
$s_R,t_R\in Hom_{C_R}(0_R,1_R)$ and a morphism $j_R\in Hom_{C_R}(1_R,0_R)$ such that $j_R\circ s_R = j_R\circ t_R = id_{0_R}$.

The category of presheaves over $C_R$ is called the category of reflexive graphs.

An object $X$ of the category $RGph$ is defined by two sets $X(0_R)$ and $X(1_R)$, two morphisms $X(s_R),X(t_R):X(1_R)\rightarrow X(0_R)$
and a morphism $X(j_R):X(0_R)\rightarrow X(1_R)$ such that $X(s_R)\circ X(j_R) = X(t_R)\circ X(j_R)=Id_{X(0_R)}$.
Let $x$ be an element of $X(0_R)$, we will often denote $X(j_R)(x)$ by $[x]$.
Geometrically, a node $x\in X(0_R)$ is represented by a point; we do not represent geometrically $X(j_R)(X(0_R))$. If  $a\in X(1_R)$ is an arc which is not an element of $X(j_R)(X(0_R))$, it is represented by a directed arrow between $X(s_R)(a)$ and $X(t_R)(a)$.

Examples of reflexive graphs are:

The reflexive dot graph $D_R$; $D_R(0_R)$ and $D_R(1_R)$ are singletons.

The reflexive arc $A_R$; $A_R(0_R)$ contains two elements $x,y$; $A_R(1_R)$ contains three elements $[x], [y]$ and $a$ such that
$A_R(j_R)(x) =[x], A_R(j_R)(y) = [y]$, $A_R(s_R)(a) = x$ and $A_R(t_R)(a) = y$.

The reflexive cycle of length $n$, $c_n^R$; $c_n^R(0_R)$ contains $n$ elements that we denote by $x^n_0,...,x^n_{n-1}$, For $i<n-1$, there is
a unique arc $a^n_i$ whose source is $x^n_i$ and whose target is $x^n_{i+1}$; there is an arc $a^n_{n-1}$ whose source is $x^n_{n-1}$ and whose target
is $x^n_0$. There exists arcs $[x^n_0],...,[x^n_{n-1}]$ such that $c_n^R(j_R)(x^n_i) = [x^n_i]$.

We are going to transport the closed models defined on $Gph$ to $RGph$.

We recall the transport theorem due to Crans, see  Cisinski [6] 1.4.23.

\medskip

{\bf Theorem 6.1.}

{\it Let,  $C$, $D$ be categories such that:

(i)  $C$ and $D$ are complete and cocomplete and $L:C\rightarrow D$ a functor which has a right adjoint $R$.

 Suppose that $C$ is endowed with a closed model structure  $(W_C,Fib_C,Cof_C)$  cofibrantly generated by $I$ and $J$ such that:
 
 (ii)  $L(I)$ and $L(J)$ allow the small element argument 
 
 (iii) for every arrow $d$ of $D$ which is the transfinite composition
of pushouts  of arrows $L(c)$ where $c$ is an element of $W_C\cap Cof_C$, the arrow $R(d)$ is a weak equivalence in $C$. 

Then
there exists a closed model structure $(W_D,Cof_D,Fib_D)$ on $D$ such that:

T1 An arrow $d$ of $D$ is in $W_D$ if and only if $R(d)$ is in $W_C$

T2 An arrow $f$ of $C$ is in $Fib_D$ if and only if $R(f)$ is in $Fib_C$.

T3 An arrow of $D$ is in $Cof_D$ if and only it has the left lifting property with respect to all elements of $W_D\cap Fib_D$.}

\medskip

 We thus deduce the following result:

\medskip

{\bf Proposition 6.1.}

{\it Let $C,D$ be categories of presheaves defined on a set. Let $(W_C,Fib_C,Cof_C)$ be a closed model defined by counting the set of objects $(X_l)_{l\in L}$ of $C$. We suppose that $Fib_C$
 is the class of all maps of $C$. Let $F:C\rightarrow D$ be a functor which has a right adjoint $G$.
 Suppose that $D$ is complete and cocomplete, then we can transfer $(W_C,Fib_C,Cof_C)$ to $D$ to obtain a closed model $(W_D,Fib_D,Cof_D)$ whose class of weak equivalences is defined by counting  the set $(F(X_l))_{l\in L}$.}

\medskip

{\bf Proof.}

The condition $(i)$ and $(ii)$  are satisfied since $C$ and $D$ are categories of presheaves defined on a set.
 Since the weak cofibrations are isomorphisms, the condition $(iii)$ is also satisfied. We deduce the class of weak equivalences of
 the closed model $(W_D,Fib_D,Cof_D)$ transfered to $D$ are morphisms $f:U\rightarrow V$ such that $G(f)$ is a weak equivalence.
  This is equivalent to saying that for every $l\in L$, the morphism of sets $Hom(X_l,G(U))\rightarrow Hom(X_l,G(V))$ which sends $h$ to $G(f)\circ h$ is an isomorphism. Since $G$ is the right adjoint of $F$, we deduce that this last condition is equivalent to saying that the morphism $Hom(F(X_l),U)\rightarrow Hom(F(X_l),V)$
  which sends $h$ to $h\circ f$ is an isomorphism. Thus $(W_D,Fib_D,Cof_D)$ is obtained by counting the family $(F(X_l))_{l\in L}$.

\medskip

We are going to apply the previous proposition to the following situation: consider the functor $f_R:C\rightarrow C_R$  defined on objects by $f_R(0) = 0_R$ and $f_R(1) = 1_R$. On morphisms, it is  defined  by
$f_R(s) = s_R$ and $f_R(t) = t_R$.

Recall that if $S$ is a presheaf defined on a category $D$, and $F:D'\rightarrow D$ a functor, the inverse image $F^*S$  of $S$ is the presheaf defined on
$D'$ such that for every object $X$ of $D'$, $F^*S(X) = S(F(X))$. When applying this construction to the functor $f_R$, we obtain that: if $X$ is a reflexive graph, $f_R^*(X)(0)=X(0_R)$ and $f_R^*(X)(1)=X(1_R)$. In particular,
 $f^*_R(D_R) = c_1$ and $f^*_R(A_R)$ is the directed graph which has two nodes $x$ and $y$, there exists an arc $a$ whose source is $x$ and whose target is $y$, there exists two loops $a_x$ such that $s(a_x)=x$ and $a_y$ such that $s(a_y)=y$.
The Proposition 5.1 p.23 of [8] insures that the functor $f^*_R$ has a left adjoint ${f_R}_*$   and a right adjoint ${f_R}_!$.

\medskip

{\bf Proposition 6.2.}

{\it The closed models of $RGph$ obtained by transferring the closed models $(W,Cof,Fib)$ which counts the cycles $(c_n)_{n\geq 0}$ and $(W',Cof',Fib')$ which counts the cycles $(c_n)_{n>0}$ to $RGph$ by the adjunction pair
$({f_R}_*,f_R^*)$ are identic.}

\medskip

{\bf Proof.}

The Proposition 6.1 implies that the transfer of $(W,Cof,Fib)$ (resp., $(W',Cof',Fib')$) on $RGph$ is the closed model defined by
counting $(c_n^R)_{n\geq 0}$ (resp., $(c_n^R)_{n\geq 1})$. Thus we have to show that a morphism $f:X\rightarrow Y$ is right orthogonal to
$i_n^R,j_n^R,n\geq 0$ if and only if it is right orthogonal to $i_n^R,j_n^R, n\geq 1$. On this purpose, it is enough to show that
if $f$ is right orthogonal to $i_n^R,j_n^R, n\geq 1$, then it is right orthogonal to $i_0^R$ and $j_0^R$. Suppose that such an $f$ is
not right orthogonal to $i_0^R$ or $j_0^R$. This equivalent to saying that $f$ does not induces a bijection between the nodes of $X$ and $Y$.
If $f(0):X(0_R)\rightarrow Y(0_R)$ is not injective, let $x,y\in X(0_R)$ such that $f(0)(x)=f(0)(y)$. There exist morphisms
$u,v:c_1^R\rightarrow X$ such that $u(0)(x^1_0) = x, v(0)(x^1_0) = y$, and $u(1)(a^1_0)=[x]$ and $v(1)(a^1_0) = [y]$. Consider the morphism
$w:c_1^R\rightarrow Y$ such that $w(0)(x^1_0) = f(0)(x)$ and $w(1)(a^1_0) = [x]$. The following
diagram does not have a filler.
$$
\begin{CD}
c_1^R+c_1^R  @>u+v>> X\\
@VV j_1^R V    @VV f V\\
c_1^R @>w>>  Y
\end{CD}
$$

This is a contradiction with the fact that $f$ is right orthogonal to $j_1^R$; thus $f(0)$ is injective.

Suppose that $f(0)$ is not surjective. Then there exists a node $y$ of $Y$ which is not in the image of $f(0)$. Let $u:c_1^R\rightarrow Y$
defined by $u(0)(x^1_0) = y$ and $u(1)(a^1_0) = [y]$. The following diagram does not have a filler:

$$
\begin{CD}
\phi  @> >> X\\
@VV i_1 V    @VV f V\\
c_1^R @>u>>  Y
\end{CD}
$$

This is in contradiction with the fact that $f$ is right orthogonal to $i_1^R$. We deduce that $f(0)$ is surjective.

\medskip

{\bf Definitions 6.1.}

Let $X$ be a reflexive graph, the cycle $f:c_n^R\rightarrow X$ is degenerated if there exists $i$ such that $f(1)(a_i^n)=[y]$
where $y$ is a node of $Y$. A cycle is nondegenerated if it is not degenerated.

\medskip

The following proposition shows that a morphism
of $W_R$ preserves the nondegenerated cycles.

\medskip

{\bf Proposition 6.3.}

{\it A weak equivalence $f:X\rightarrow Y$ of $RGph$  induces a bijection  on nondegenerated cycles.}

\medskip

{\bf Proof.}

Suppose that the image of a cycle $u:c_n^R\rightarrow X$ is degenerated. This implies that there exists a cycle $v:c_{n-1}^R\rightarrow Y$ which has the same image than $u$ and such that there exists a commutative  diagram:

$$
\begin{CD}
\phi  @> >> X\\
@VV  V    @VV f V\\
c_{n-1}^R @>v>>  Y
\end{CD}
$$

which has a filler $w:c_{n-1}^R\rightarrow X$, and there exists a degenerated morphism
$h:c_n^R\rightarrow c^R_{n-1}$ such that $ f\circ w\circ h =f\circ u$.
Since the image of $w$ and the image of $u$ are different, we deduce that $f$ does not induces an injection on $n$-cycles. This is a contradiction with the fact that $f$ is a weak equivalence.

\medskip

There exist morphisms which induces bijection on nondegenerated cycles, but which are not weak equivalences an example is the canonical morphism $f:A_R\rightarrow D_R$. The following result can be compared to [10] Theorem 4.9:

\medskip

{\bf Proposition 6.4.}

{\it Let $W'_R$ be the class of morphisms of $RGph$ which induce a bijection on nondegenerated cycles. There does not exist
a closed model whose class of weak equivalences is $W'_R$.}

\medskip

{\bf Proof.}

Suppose that such a closed model exists. Consider the canonical morphism $f:A_R\rightarrow D_R$, we can write $f=g\circ h$ where $g$ is a weak fibration and $h$ a cofibration, the $2$-$3$ property implies that $h$ is a weak cofibration.
Write $g:X\rightarrow D_R$, suppose that the cardinality of $X(0_R)$ is superior or equal to $2$. Let $l:c_1^R\rightarrow D_R$, the pullback of $l$ by $g$ is not a weak equivalence since its domain contains at least two distinct subgraphs isomorphic to $c_1^R$, this implies that the cardinal of $X(0_R)$ is $1$ and henceforth the cardinal of $X(1_R)$ is $1$ since $g$ is a weak equivalence; thus $g$ is the identity. We deduce that $f=h$ is a weak cofibration.

Let $Y$ be the reflexive graph such that $Y(0_R)$ contains two elements $u$ and $v$, $Y(1_R)$ contains $[u], [v]$ and two elements $c,d$ such that $X(s_R)(c)=X(s_R)(d) =u$ and $X(t_R)(c) = X(t_R)(d) = v$. Consider the morphism $k:A_R\rightarrow Y$ such that $k(0)(x)=u, k(0)(y) = v$ and $k(1)(a) =c$.
The image of the pushout $m$ of $f$ by $k$ is $c_1^R$. This implies that $m$ is not weak equivalence. This is a contradiction with the fact that the pushout of a weak cofibration is a weak cofibration.

\medskip

We will show now that some properties of the closed model defined on $RGph$ similar to the properties of the closed model $(W,Fib,Cof)$ defined on $Gph$. A reflexive graph $X$ is strongly connected if and only if for every nodes $x$ and $y$ of $X$, there exists a reflexive cycle $f:c_n^R\rightarrow X$ such that the image of $f(0)$ contains $x$ and $y$.

\bigskip

{\bf Proposition 6.5.}

{\it A strongly connected reflexive graph $X$ is cofibrant.}

\medskip

{\bf Proof.}

Let $X$ be a strongly connected reflexive graph $X$. There exists a graph $X'$ in $Gph$ such that ${f_R}_*(X')=X$; $X'(0) = X(0_R)$
and $X'(1)$ is $X(1_R)-\{ [x], x\in X(0)\}$. The graph $X'$ is also strongly connected, thus it is a cofibrant object
of $(W,Cof,Fib)$. Since the map $c_{X'}:\phi\rightarrow X'$ is a cofibration, this implies that $c_{X'}$ is an element of $cell(i_n,j_n,n\geq 0)$. We deduce that $c_X:\phi\rightarrow X$ is an element of $cell(i_n^R,j_n^R,n\geq 0)$ since $X={f_R}_*(X')$ and left adjoint preserve colimits and henceforth that $X$ is a cofibrant object. 

\medskip

{\bf Proposition 6.6.}

{\it A morphism $f:X\rightarrow Y$ between two reflexive graphs is a weak equivalence if and only if it induces a bijection between strongly connected components and its restriction to each strongly connected component is an isomorphism onto a strongly connected component of $Y$.}

\medskip

{\bf Proof.}

Let $f:X\rightarrow Y$ be a weak equivalence of the closed model defined on $RGph$. The morphism ${f_R}^*(f)$ is also a weak equivalence. The Theorem 4.1 implies that it induces a bijection between the strongly connected components of ${f_R}^*(X)$ and ${f_R}^*(Y)$ and the restriction of ${f_R}^*(f)$ to each connected component of ${f_R}^*(X)$ is an isomorphism. Remark that ${f_R}^*(X)(0)=X(0_R)$ and ${f_R}^*(X)(1)=X(1_R)$, since ${f_R}^*$ is just the forgetful functor. This implies that the strongly connected components of ${f_R}^*(X)$ are of the form $V={f_R}^*(U)$ where $U$ is a strongly connected component of $X$ and that $f$ induces a bijection between strongly connected components and its restriction to each strongly connected component is an isomorphism onto a strongly connected component of $Y$.

\medskip

\bigskip

{\bf References.}

\bigskip

[1] Beke, T. (2000). Sheafifiable homotopy model categories. In Mathematical Proceedings of the Cambridge Philosophical Society (Vol. 129, No. 03, pp. 447-475). Cambridge University Press.
 
 \smallskip

[2] Bisson, T.,  Tsemo, A. (2009). A homotopical algebra of graphs related to zeta series. Homology, Homotopy and Applications, 11(1), 171-184.

\smallskip

[3] Bisson, T.,  Tsemo, A. (2011). Symbolic dynamics and the category of graphs. Theory and Applications of Categories, 25(22), 614-640.

\smallskip

[4] Bisson, T.,  Tsemo, A. (2011). Homotopy equivalence of isospectral graphs. New York J. Math, 17, 295-320.

\smallskip

[5] Brin, S.,  Page, L. (2012). Reprint of: The anatomy of a large-scale hypertextual web search engine. Computer networks, 56(18), 3825-3833.
\smallskip

[6] Cisinski D.C.,(2006)  Les pr\'efaisceaux comme type d'homotopie, Ast\'erisque, Volume 308, Soc. Math. France.

\smallskip

[7] Euler, L. (1741). Solutio problematis ad geometriam situs pertinentis. Commentarii academiae scientiarum Petropolitanae, 8, 128-140.

\smallskip

[8] Artin, M., Grothendieck, A.,  Verdier, J. L. (1972). Th\'eorie des topos et cohomologie \'etale des sch\'emas. Tome 1. Lecture notes in mathematics, 269.

\smallskip

[9] Hirschhorn, P. S. (2009). Model categories and their localizations (No. 99). American Mathematical Soc.

\smallskip

[10] Tsemo, A. (2013). Applications of closed models defined by counting to graph theory and topology. arXiv preprint arXiv:1308.3983.

\end{document}